\documentclass{article}
\usepackage[utf8]{inputenc}
\usepackage{cite}
\usepackage{url}
\usepackage{caption, subcaption}
\usepackage{hyperref}
\hypersetup{
    colorlinks=true,
    linkcolor=red,
    filecolor=blue,      
    urlcolor=blue,
    citecolor=blue}
\usepackage{amsmath,amssymb,amsfonts}
\usepackage[noend]{algorithmic}
\usepackage{algorithm}
\makeatletter
\newcommand{\removelatexerror}{\let\@latex@error\@gobble}
\makeatother
\usepackage{multirow}
\usepackage{graphicx}
\usepackage{textcomp}
\usepackage{enumitem}
\usepackage{xcolor}
\usepackage{color}
\usepackage{bm}
\usepackage{bbm}
\usepackage{todonotes}
\newcommand{\x}{\bm{x}}
\newcommand{\y}{\bm{y}}
\newtheorem{remark}{Remark}
\newtheorem{definition}{Definition}

\newtheorem{example}{Example}

\providecommand{\keywords}[1]
{
  \small	
  \textbf{Key words.} #1
}

\providecommand{\MSC}[1]
{
  \small	
  \textbf{MSC codes.} #1
}

\title{A Linear Complexity Algorithm for Optimal Transport Problem with Log-type Cost}
\author{Ziyuan Lyu\footnote{Department of Mathematical Sciences, Tsinghua University, Beijing 100084, China
(lvzy21@mails.tsinghua.edu.cn).}, Zihao Wang\footnote{Corresponding author. Department of Computer Science and Engineering, Hong Kong University of Science and Technology, Clear Water Bay, Hong Kong SAR, China (zwanggc@cse.ust.hk).}, Hao Wu\footnote{Department of Mathematical Sciences, Tsinghua University, Beijing 100084,
China (hwu@tsinghua.edu.cn).}, and Shuai Yang\footnote{Department of Mathematical Sciences, Tsinghua University, Beijing 100084, China
(s-yang21@mails.tsinghua.edu.cn).}}

\begin{document}

\maketitle
\begin{abstract}
    In [Q. Liao et al., Commun. Math. Sci., 20(2022)], a linear-time Sinkhorn algorithm is developed based on dynamic programming, which significantly reduces the computational complexity involved in solving optimal transport problems. However, this algorithm is specifically designed for the Wasserstein-1 metric. We are curious whether the preceding dynamic programming framework can be extended to tackle optimal transport problems with different transport costs. Notably, two special kinds of optimal transport problems, the Sinkhorn ranking and the far-field reflector and refractor problems, are closely associated with the log-type transport costs. Interestingly, by employing series rearrangement and dynamic programming techniques, it is feasible to perform the matrix-vector multiplication within the Sinkhorn iteration in linear time for this type of cost. This paper provides a detailed exposition of its implementation and applications, with numerical simulations demonstrating the effectiveness and efficiency of our methods.
\end{abstract}

\keywords{Optimal Transport, Sinkhorn Algorithm, Log-type Transport Cost, Sinkhorn Ranking, Far-field Reflector and Refractor}

\MSC{49M25; 65K10}

\section{Introduction}
In this study, we aim to design a linear complexity Sinkhorn-type algorithm for the optimal transport problem with log-type transport cost. The Sinkhorn algorithm is a popular numerical algorithm for solving optimal transport problems\cite{cuturi2013sinkhorn, sinkhorn1967diagonal} and has become one of the most competitive methods among numerous optimal transport numerical algorithms\cite{li2020asymptotically, bertsekas1988auction, glimm2013iterative, benamou2000computational, lee2021generalized, benamou2014numerical, prins2015least, hu2022efficient} due to its simple iterative form and its applicability to various transport costs. As is well known, the Sinkhorn algorithm alternately updates scaling variables, which includes repeated multiplication of the kernel matrix and the scaling variables and leads to $O(N^2)$ time complexity, where $N$ is the number of discrete points. In practical scenarios, $N$ often greatly exceeds $10^4$, thereby limiting the applicability of the Sinkhorn algorithm. In\cite{liao2022fast} and\cite{liao2022fast2}, dynamic programming techniques are applied to develop methods to implement the Sinkhorn algorithm in $O(N)$ time by using the special structure of the kernel matrix derived from the Wasserstein-1 metric. On this basis, our question is whether valuable matrix structure information can be extracted for optimal transport problems with other types of transport costs, with the goal of accelerating the Sinkhorn algorithm and enhancing its application to related problems.

Optimal transport problem was originally introduced by Monge in 1781\cite{monge1781memoire}. He considered a practical problem of how to transport the earth from a given place to another given place with the least total cost, where the transport cost is $c(\x,\y)=\|\x-\y\|$. However, Monge's problem is difficult to solve, as it is challenging to answer whether a minimizer of the problem exists. Thanks to Kantorovich's generalization to Monge's problem in 1942\cite{kantorovich1942transfer}, manageable solutions can be found and studied. In Kantorovich's problem, the transport cost is extended to more general cases\cite{santambrogio2015optimal}, with some typical examples including
\begin{enumerate}[label={(\arabic*)}]
    \item $c(\x,\y)=d(\x,\y)^p$ for $p \ge 1$, where $d(\cdot,\cdot)$ is a distance on a Polish metric space. This cost is used to define the Wasserstein-$p$ metric\cite{villani2009optimal}. The Wasserstein-$p$ metric has many applications in various fields. For example, the Wasserstein-1 metric is applied in integrated information theory\cite{oizumi2014phenomenology}; the Wasserstein-2 metric has a nice displacement interpolation property\cite{mccann1997convexity} and is widely used in computer graphics\cite{bonneel2011displacement, solomon2015convolutional}. 
    \item $c(\x,\y)=-\log(1-\langle \x,\y\rangle)$ and $c(\x,\y)=-\log(1-\kappa\langle \x,\y\rangle)$. These are the transport costs that appear in the far-field reflector problem and the far-field refractor problem, where $0<\kappa < 1$ is the relative refractive index \cite{wang2004design, gutierrez2009refractor}. 
    \item $c(\x_1,\x_2,\cdots, \x_N)=\sum_{1\le i<j\le N}1/\|\x_i-\x_j\|$. This is the Coulomb cost, which arises in many-electron quantum mechanics. It is related to the $N$-body optimal transport problem\cite{cotar2013density, cotar2015infinite}. 
\end{enumerate}
Currently, for the Wasserstein-1 and Wasserstein-2 metrics, there has been extensive research on numerical methods and acceleration techniques\cite{li2018parallel, korotin2019wasserstein, mallasto2017learning, arjovsky2017wasserstein, solomon2015convolutional, can2019accelerated, liu2021multilevel, liao2022fast, liao2022fast2}. However, research on numerical methods for other types of transport costs is somewhat limited.

To our best knowledge, the optimal transport problem with log-type transport cost was first introduced in\cite{wang2001monge, wang2004design}, where the far-field reflector problem is proved to be an optimal transport problem with transport cost $c(\x,\y)=-\log(1-\langle \x,\y\rangle)$. Following this work, Guti{\'e}rrez et al. proved that the far-field refractor problem is also an optimal transport problem with similar transport cost\cite{gutierrez2009refractor}. Later on more complex optical systems have been studied, including the system of a reflector pair\cite{glimm2010rigorous} and the system of a refractor pair\cite{oliker2018beam}. For this kind of problem, there are linear assignment\cite{doskolovich2019optimal} and linear programming methods\cite{glimm2003optical, wang2004design, oliker2018beam}. Besides, the Sinkhorn algorithm was considered to solve the reflector problem\cite{benamou2020entropic}, and its convergence as $N$ and the inverse of the regularization parameter $1/\varepsilon$ jointly tend to infinity was proved in\cite{berman2020sinkhorn}. There are also some newly developing PDE based approaches that solve the reflector problem through its corresponding Monge-Amp{\`e}re type equation, including a generalized least-squares method\cite{romijn2020inverse} and a convergent finite difference scheme\cite{hamfeldt2021convergent}. In recent years, some new forms of log-type transport costs have emerged, such as the log-determinant cost\cite{altschuler2021hardness}. The transport cost of the form $c(\x, \y) = -\varepsilon\log \langle \varphi(\x), \varphi(\y)\rangle$ was proposed and applied to approximate the entropy regularized optimal transport distance, where $\varphi$ is a positive feature map. With this log-type transport cost, the time complexity of each Sinkhorn iteration can be reduced to $O(sN)$, where $s$ denotes the dimension of the feature space\cite{scetbon2020linear}. Additionally, the optimal transport problem with log-type transport cost is also closely related to differentiable ranking\cite{cuturi2019differentiable}.

In this paper, we consider the log-type transport cost that is of the form $c(\x,\y)=-\log(P(\x,\y))$, where $P$ is a polynomial of the entries of $\x,\y$. The transport costs derived from the far-field reflector problem and the far-field refractor problem with the relative refractive index $0<\kappa < 1$\cite{gutierrez2009refractor} are both in this form. Some transport costs in this form can be used for differentiable ranking. When using the Sinkhorn algorithm to solve an optimal transport problem with this kind of transport cost, each element of the corresponding kernel matrix can be written as a polynomial expansion under the assumption that the regularization parameter $\varepsilon=1/L$ is the inverse of a positive integer $L$.\footnote{The solution to the entropy regularized optimal transport problem converges to a solution to the unregularized optimal transport problem as $\varepsilon\rightarrow 0$. In actual use, one is free to set the value of $\varepsilon$.} In this way, series rearrangement and dynamic programming techniques can be applied to accelerate the matrix-vector multiplication within the Sinkhorn iteration. From the perspective of time complexity, for the $d$-dimensional case, our method performs each Sinkhorn iteration with time complexity $O((ML)^dN)$, where $M$ is the maximum of the exponents of the entries of $\x,\y$ among all the monomials of $P$. When $(ML)^d \ll N$, which naturally holds true since $N$ is commonly significantly large in the $d$-dimensional case (see Remark~\ref{rmk::M=1or2} for a detailed discussion), the new method exhibits a substantial efficiency advantage over the original Sinkhorn algorithm. We name our algorithm Fast Sinkhorn for optimal transport with Log-type cost (FSL).

It is worth mentioning that, our work is motivated by the fast Sinkhorn algorithms\cite{liao2022fast, liao2022fast2} and fast matrix-vector multiplication algorithms\cite{koev1999matrices, indyk2022faster}. For fast Sinkhorn algorithms~\cite{liao2022fast, liao2022fast2}, the $N\times N$ kernel matrix is required to be a lower/upper-collinear triangular matrix, which requires an uniform discretization. However, the kernel matrix we consider in this paper applies to general discretization, which follows the recent advance in general distance matrix~\cite{indyk2022faster}, or in general, the matrix with displacement structure~\cite{koev1999matrices}. Compared to those works only targeting on matrix-vector multiplication~\cite{indyk2022faster,koev1999matrices}, we systematically leverage the Sinkhorn algorithm and incorporate specialized techniques designed for the kernel matrix associated with log-type transport cost to achieve linear-time complexity.

Our FSL algorithm provides meaningful practical utility and exhibits promising application potential. As an example, we give a new approach to compute the Sinkhorn ranking operator\cite{cuturi2019differentiable} (Section~\ref{subsec::1d_FS4Implementation}), wherein a novel log-type cost function is adopted to define this operator, and the FSL algorithm is applied to greatly reduce the computational cost. The Sinkhorn ranking operator belongs to the family of differentiable ranking operators, serving as a differentiable proxy to the usual ranking operator. While ranking is important in machine learning, the non-differentiability of this operation restricts its usage in end-to-end learning\cite{cuturi2019differentiable}. Differentiable proxies facilitate the end-to-end training using ranking, with many applications including learning-to-rank models\cite{swezey2021pirank} and neural network-based $k$-nearest neighbor classifiers\cite{cuturi2019differentiable, xie2020differentiable}. Researchers have proposed a variety of differentiable ranking operators, for instance, the SoftRank approach that using random perturbation technique\cite{taylor2008softrank}, and the method that utilizing the pairwise difference matrix\cite{qin2010general}. Among them, the Sinkhorn ranking operator\cite{cuturi2019differentiable} has gained increasing attention due to its broad applicability in recent years. It is well known that usual ranking costs $O(N\log N)$ time. To our best knowledge, the best-known result regarding the time complexity of differentiable ranking operator is also $O(N \log N)$\cite{blondel2020fast}. As for the computational complexity of the original Sinkhorn ranking operator, it takes $O(N^2)$ time to perform each iteration when both the number of values to be ranked and the number of auxiliary values are $N$\cite{cuturi2019differentiable}. Our approach, however, can achieve linear-time complexity for each iteration. Numerical simulations illustrate the remarkable efficiency of our approach in computing the Sinkhorn ranking operator (Section~\ref{subsec::1dexperiments}). Furthermore, the results show that the Sinkhorn ranking operator defined by our log-type cost function can give better performance than the one defined by the squared distance cost function as in\cite{cuturi2019differentiable} in some scenarios. Consequently, our approach is potentially a better choice for obtaining such operator.

We also consider the application of the FSL algorithm to the reflector cost, the log-type transport cost that derived from the far-field reflector problem\cite{wang2004design}; and to the refractor cost, the log-type transport cost that comes from the far-field refractor problem with the relative refractive index $0<\kappa < 1$\cite{gutierrez2009refractor} (Section~\ref{sec::app_2d}). Each of these problems wants to find a surface that redirects the light from the point source to the target region in a far-field, and obtain the prescribed illumination intensity on that region. One can recover their generalized solutions by solving the corresponding dual Kantorovich's problems. When using the Sinkhorn algorithm to solve the corresponding entropy regularized optimal transport problem, the computational cost of each iteration is $O(N^2)$, which is demanding as $N$ is often quite large in these problems\cite{benamou2020entropic}. In order to alleviate this issue, we utilize our FSL algorithm to develop a linear-time algorithm for these two transport costs. 

The rest of this paper is organized as follows. In Section~\ref{sec::FS4}, we briefly introduce the Sinkhorn algorithm and then present our FSL algorithm. In Section~\ref{sec::app_1d}, we introduce how to define the Sinkhorn ranking operator by our proposed log-type cost function and how to implement the FSL algorithm to speed up the computation, and present some numerical simulations of this approach. In Section~\ref{sec::app_2d}, we apply the FSL algorithm to the reflector and refractor costs. Finally, the conclusions are given in Section~\ref{sec::conclusion}.

\section{Fast Sinkhorn for Log-type Transport Cost}\label{sec::FS4}
In this section, we first give some notational conventions and the definition of the log-type transport cost. After that, we briefly recall the Sinkhorn algorithm. Then, we detail our FSL algorithm.

\subsection{Notations and Definitions}
Throughout this paper, we consider the optimal transport problem on two discrete finite spaces $X,Y \subset \mathbb{R}^d$, say $X=\{\x_1,\x_2,\cdots,\x_N\}$, $Y=\{\y_1,\y_2,\cdots,\y_N\}$, and $\x_i = (x_{i,1},x_{i,2},\cdots,x_{i,d})^\top$, $\y_i = (y_{i,1},y_{i,2},\cdots,y_{i,d})^\top$, $i=1,\cdots,N$. 

We use a vector to represent a function or measure on $X$, and use a matrix to represent a function or measure on $X \times Y$. For example, a function $\bm{f}:X \rightarrow \mathbb{R}$ is denoted by $\bm{f}=(f_1,f_2,\cdots,f_N)^\top \in \mathbb{R}^N$, where $f_i = \bm{f}(\x_i)$; a function $\bm{G}:X\times Y \rightarrow \mathbb{R}$ is denoted by $\bm{G}=(G_{ij}) \in \mathbb{R}^{N\times N}$, where $G_{ij}=\bm{G}(\x_i,\y_j)$. 

The set of probability measures on $X$ is denoted by $\mathcal{P}(X)$. We use $\langle \cdot,\cdot\rangle$ to represent the standard Euclidean inner product. The vector of ones in $\mathbb{R}^N$ is written as $\bm{1}_N = (1,1,\cdots,1)^\top \in \mathbb{R}^N$. For two vectors $\bm{f},\bm{g} \in \mathbb{R}^N$ (suppose $\bm{g}$ is a positive vector), the operator $\oslash$ represents pointwise division, e.g., $\bm{f} \oslash \bm{g} \in \mathbb{R}^N$, $(\bm{f} \oslash \bm{g})_i = f_i/g_i$. The sets of positive and non-negative real numbers are denoted by $\mathbb{R}_{++}$ and $\mathbb{R}_+$, respectively.

For a positive integer $k$ with $k=k_1+k_2+\cdots+k_m$, where $m \ge 2$ and $k_1,\cdots,k_m \ge 0$ are integers, the multinomial coefficient is given by
\begin{align*}
    \binom{k}{k_1\,\,k_2\,\cdots\, k_m} = \frac{k!}{k_1!k_2!\cdots k_m!}.
\end{align*}

As mentioned before, the log-type transport cost that we concern is of the form $c(\x,\y)=-\log(P(\x,\y))$. More precisely, its definition is as follows: 
\begin{definition}\label{def::log-type cost}
    Given a polynomial $P$ of $2d$ variables with coefficients in $\mathbb{R}$ that satisfying
\begin{align}
    P(x_{i,1},\cdots,x_{i,d},y_{j,1},\cdots,y_{j,d})\in (0,1), \quad\forall i,j=1,2,\cdots,N,
    \label{eq::P_assume}
\end{align}
the log-type transport cost $\bm{C} = (C_{ij})$ associated with $P$ is given by
\begin{align*}
    C_{ij} = -\log(P(x_{i,1},\cdots,x_{i,d},y_{j,1},\cdots,y_{j,d})),\quad \forall i,j=1,2,\cdots,N.
\end{align*}
The polynomial $P$ has the form 
\begin{multline}
    P(x_{i,1},\cdots,x_{i,d},y_{j,1},\cdots,y_{j,d})=\sum_{\zeta_1 ,\cdots,\zeta_d ,\nu_1 ,\cdots,\nu_d =0}^M a_{\zeta_1 ,\cdots,\zeta_d ,\nu_1 ,\cdots,\nu_d }\prod_{k=1}^d x_{i,k}^{\zeta_k }y_{j,k}^{\nu_k },
    \label{eq::P}
\end{multline}
where $a_{\zeta_1 ,\cdots,\zeta_d ,\nu_1 ,\cdots,\nu_d }$ is the coefficient of the monomial $\prod_{k=1}^d x_{i,k}^{\zeta_k }y_{j,k}^{\nu_k }$, and $M$ is the maximal positive integer such that $\exists \,a_{\zeta_1 ,\cdots,\zeta_d ,\nu_1 ,\cdots,\nu_d } \ne 0$ with $\zeta_1 ,\cdots,\zeta_d ,\nu_1 ,\cdots,\nu_d  \le M$ and $\prod_{k=1}^d(M-\zeta_k )(M-\nu_k )=0$.
\end{definition}

\subsection{Sinkhorn Algorithm}
For two probability measures $\bm{a} \in \mathcal{P}(X)$ and $\bm{b} \in \mathcal{P}(Y)$ and a cost matrix $\bm{C} = (C_{ij}) \in \mathbb{R}^{N\times N}$, the Kantorovich's optimal transport problem between $\bm{a}$ and $\bm{b}$ is the following optimization problem 
\begin{align}
    I(\bm{a},\bm{b}) &\triangleq \min_{\bm{\Gamma} \in \Pi(\bm{a},\bm{b})} \langle \bm{\Gamma},\bm{C}\rangle,\nonumber\\
    \Pi(\bm{a},\bm{b}) &\triangleq \{\bm{\Gamma} \in \mathcal{P}(X \times Y)\,|\,\bm{\Gamma} \bm{1}_N=\bm{a}, \bm{\Gamma}^\top \bm{1}_N=\bm{b}\},
 \label{eq::KP}
\end{align}
where $\Pi(\bm{a},\bm{b})$ is the set of transport plans between $\bm{a}$ and $\bm{b}$. It can be seen that $\bm{a}\bm{b}^\top \in \Pi(\bm{a},\bm{b})$, hence $\Pi(\bm{a},\bm{b})$ is non-empty. Since $\langle \bm{\Gamma},\bm{C}\rangle$ is continuous on $\Pi(\bm{a},\bm{b})$, and that $\Pi(\bm{a},\bm{b})$ is compact, the Kantorovich's problem~\eqref{eq::KP} has at least one optimal solution.

To solve~\eqref{eq::KP}, the Sinkhorn algorithm\cite{cuturi2013sinkhorn, sinkhorn1967diagonal} considers the following entropy regularized optimal transport problem
\begin{align}
    \min_{\bm{\Gamma} \in \Pi(\bm{a},\bm{b})} \langle \bm{\Gamma},\bm{C}\rangle + \varepsilon \sum_{i,j=1}^N \Gamma_{ij}\log\left(\Gamma_{ij}\right),
 \label{eq::entropyKP}
\end{align}
where $\varepsilon > 0$ is called the regularization parameter. Starting from some positive function $\bm{\phi}^{(0)}$ on $X$, the Sinkhorn algorithm iteratively computes the vectors
\begin{align}
        \bm{\psi}^{(t+1)} = \bm{b} \oslash \left(\bm{K}^\top \bm{\phi}^{(t)}\right), \quad \bm{\phi}^{(t+1)} = \bm{a} \oslash \left(\bm{K} \bm{\psi}^{(t+1)}\right),
        \label{eq::SinkhornIterations}
\end{align}
where $t$ is the iteration counter, $\bm{K}\triangleq \mbox{exp}(-\bm{C}/\varepsilon)$ is the kernel associated with $\bm{C}$ and $\varepsilon$, i.e., 
\begin{align*}
    K_{ij} = e^{-\frac{C_{ij}}{\varepsilon}},\quad i,j=1,2,\cdots,N.
\end{align*}
When the termination criterion is satisfied, the vectors $\bm{\phi}^{(t)}$ and $\bm{\psi}^{(t)}$ are returned. From them, we can obtain the approximate optimal transport plan $\mbox{diag}(\bm{\phi}^{(t)}) \bm{K} \mbox{diag}(\bm{\psi}^{(t)})$, and use $\langle \mbox{diag}(\bm{\phi}^{(t)}) \bm{K} \mbox{diag}(\bm{\psi}^{(t)}), \bm{C} \rangle$ to approximate $I(\bm{a},\bm{b})$. 

The main computational cost of each Sinkhorn iteration comes from the two matrix-vector multiplications, which takes $O(N^2)$. The pseudo-code of the Sinkhorn algorithm is presented in Algorithm~\ref{alg::Sinkhorn}.

\begin{figure}[!t]
  \renewcommand{\algorithmicrequire}{\textbf{Input:}}
  \renewcommand{\algorithmicensure}{\textbf{Output:}}
  \removelatexerror
  \begin{algorithm}[H]
    \caption{Sinkhorn Algorithm}
    \label{alg::Sinkhorn}
    \begin{algorithmic}[1]
      \REQUIRE $\varepsilon, \bm{a}, \bm{b}, \bm{\phi}^{(0)}, \bm{C}$
      \STATE $\bm{K}\leftarrow  e^{-\bm{C}/\varepsilon}$; $t \leftarrow 0$
      \WHILE{not terminate}
        \STATE $\bm{\psi}^{(t+1)} \leftarrow \bm{b} \oslash \left(\bm{K}^\top \bm{\phi}^{(t)}\right)$
        \STATE $\bm{\phi}^{(t+1)} \leftarrow \bm{a} \oslash \left(\bm{K} \bm{\psi}^{(t+1)}\right)$
        \STATE $t\leftarrow t+1$
      \ENDWHILE
      \RETURN $\bm{\phi}^{(t)}, \bm{\psi}^{(t)}$
    \end{algorithmic}
  \end{algorithm}
\end{figure}

\subsection{FSL Algorithm}\label{subsec::FSL}
We now present our FSL algorithm for reducing the computational cost of each Sinkhorn iteration with respect to the log-type transport cost. 

For simplicity of presentation, we will restrict ourselves to the one-dimensional case, i.e., $d=1$, and denote by $x_j=x_{j,1}$ and $y_j=y_{j,1}$. We further assume that $L\triangleq 1/\varepsilon \in \mathbb{N}$. In this case, the elements of the kernel $\bm{K}$ can be written as polynomial expansions
\begin{multline}
    K_{ij} = P(x_i,y_j)^{\frac{1}{\varepsilon}} 
    = \left(\sum_{\zeta=0}^M\sum_{\nu=0}^M a_{\zeta,\nu}x_i^{\zeta}y_j^{\nu}\right)^{L}\\
    =\sum_{\zeta=0}^{ML}\sum_{\nu=0}^{ML} b_{\zeta,\nu}x_i^{\zeta}y_j^{\nu},\quad i,j=1,2,\cdots,N,
    \label{eq::b}
\end{multline}
where $b_{\zeta,\nu}$ is the coefficient of $x_i^{\zeta}y_j^{\nu}$ in $K_{ij}$. For a vector $\bm{\xi} \in \mathbb{R}^N$, each entry of $\bm{K} \bm{\xi}$ can be computed in the following way:
\begin{multline}
        (\bm{K} \bm{\xi})_i =\sum_{j=1}^N K_{ij} \xi_j
        =\sum_{j=1}^N \left(\sum_{\zeta=0}^{ML}\sum_{\nu=0}^{ML} b_{\zeta,\nu}x_i^{\zeta}y_j^{\nu}\right)\xi_j\\
        =\sum_{\zeta=0}^{ML}x_i^{\zeta}\sum_{\nu=0}^{ML} b_{\zeta,\nu}  \left(\sum_{j=1}^N y_j^{\nu}\xi_j\right), \quad i=1,2,\cdots, N.
        \label{eq::FSL_1}
\end{multline}
For convenience, let us define
\begin{align}
    S_{\zeta} \triangleq \sum_{\nu=0}^{ML} b_{\zeta,\nu}  \left(\sum_{j=1}^N y_j^{\nu}\xi_j\right), \quad \zeta=0,1,\cdots,ML.
    \label{eq::S}
\end{align}
And for $z=x$ or $y$, we introduce
\begin{align}
    \bm{A}^z &= (A^z_{\zeta j})\triangleq \left(
    \begin{array}{cccc}
        1 & 1 & \cdots & 1\\
        z_1 & z_2 & \cdots & z_N\\
        \vdots & \vdots & \ddots & \vdots\\
        (z_1)^{ML} & (z_2)^{ML} & \cdots & (z_N)^{ML}
    \end{array}
    \right)\in \mathbb{R}^{(ML+1)\times N},\nonumber\\
    A^z_{\zeta j} &\triangleq (z_j)^\zeta,\quad \zeta=0,1,\cdots,ML,\,\, j=1,2,\cdots,N.
    \label{eq::A}
\end{align}
Substituting~\eqref{eq::A} into~\eqref{eq::S}, and plugging~\eqref{eq::S} into~\eqref{eq::FSL_1}, we obtain that
\begin{align}
    S_\zeta = \sum_{\nu=0}^{ML} b_{\zeta,\nu}  \left(\bm{A}^y\bm{\xi}\right)_{\nu}, \quad \zeta=0,1,\cdots,ML,
    \label{eq::Sshort}
\end{align}
and
\begin{align}
    (\bm{K} \bm{\xi})_i = \sum_{\zeta=0}^{ML}A^x_{\zeta i}S_{\zeta}, \quad i=1,2,\cdots, N.
    \label{eq::FSL_2}
\end{align}
Based on~\eqref{eq::Sshort} and~\eqref{eq::FSL_2}, instead of computing $\bm{K} \bm{\xi}$ directly, FSL obtains it in the following three steps:
\begin{enumerate}
    \item Compute $\bm{A}^y\bm{\xi}$, which takes $O(MLN)$.
    \item Compute $S_\zeta$ for $\zeta=0,1,\cdots,ML$ by substituting the obtained vector from the first step for $\bm{A}^y\bm{\xi}$ in~\eqref{eq::Sshort}, which takes $O(M^2L^2)$.
    \item Compute each entry of $\bm{K} \bm{\xi}$ using the obtained values $S_\zeta$ and~\eqref{eq::FSL_2}, which takes $O(MLN)$.
\end{enumerate}

It should be noted that, FSL constructs the matrices $\bm{A}^x$ and $\bm{A}^y$ and evaluates the coefficients $b_{\zeta,\nu}$ before the main iterations~\eqref{eq::SinkhornIterations}. In this way, we can exploit the existing values of $\bm{A}\bm{\xi}$ and $S_\zeta$, so as to avoid the re-computation of these intermediate results. The vector $\bm{K}^\top \bm{\xi}$ can be obtained similarly, and we will not repeat here. 

It can be found that, the computational cost of updating $\bm{\phi}^{(t+1)}$ or $\bm{\psi}^{(t+1)}$ in FSL is reduced from $O(N^2)$ to $O(MLN)$ when $ML \ll N$. The pseudo-code of the FSL algorithm for the one-dimensional case is presented in Algorithm~\ref{alg::FSL_1d}.
\begin{figure}[!t]
  \renewcommand{\algorithmicrequire}{\textbf{Input:}}
  \renewcommand{\algorithmicensure}{\textbf{Output:}}
  \removelatexerror
  \begin{algorithm}[H]
    \caption{1D FSL Algorithm }
    \label{alg::FSL_1d}
    \begin{algorithmic}[1]
      \REQUIRE $L = \frac{1}{\varepsilon}, \bm{a},\bm{b}, \bm{\phi}^{(0)}, P(\cdot,\cdot)$
      \STATE $t \leftarrow 0$; $\bm{r} \leftarrow \bm{0}_N$; $\bm{\eta} \leftarrow \bm{0}_{ML+1}$
      \STATE Construct $\bm{A}^x$ and $\bm{A}^y$ by~\eqref{eq::A}; compute the coefficients $b_{\zeta,\nu}$ in~\eqref{eq::b}
      \WHILE{not terminate}
        \STATE $\bm{\eta} \leftarrow \bm{A}^x\bm{\phi}^{(t)}$
        \FOR{$\nu=0$ to $ML$}
            \STATE $S_\nu \leftarrow \sum_{\zeta=0}^{ML} b_{\zeta,\nu}  \eta_{\zeta}$
        \ENDFOR
        \FOR{$i=1$ to $N$}
            \STATE $r_i \leftarrow \sum_{\nu=0}^{ML}A^y_{\nu i}S_{\nu}$
        \ENDFOR
        \STATE $\bm{\psi}^{(t+1)} \leftarrow \bm{b} \oslash \bm{r}$
        
        \STATE $\bm{\eta} \leftarrow \bm{A}^y\bm{\psi}^{(t+1)}$
        \FOR{$\zeta=0$ to $ML$}
            \STATE $S_\zeta \leftarrow \sum_{\nu=0}^{ML} b_{\zeta,\nu}  \eta_{\nu}$
        \ENDFOR
        \FOR{$i=1$ to $N$}
            \STATE $r_i \leftarrow \sum_{\zeta=0}^{ML}A^x_{\zeta i}S_{\zeta}$
        \ENDFOR
        \STATE $\bm{\phi}^{(t+1)} \leftarrow \bm{a} \oslash \bm{r}$
        \STATE $t\leftarrow t+1$
      \ENDWHILE
      \RETURN $\bm{\phi}^{(t)}, \bm{\psi}^{(t)}$
    \end{algorithmic}
  \end{algorithm}
\end{figure}

\begin{remark}\label{rmk::M=1or2}
    In practice, the maximal degree $M$ of the polynomial $P$ in one variable is often as low as $1$ or $2$. Recall that $L = 1/\varepsilon$. In the numerical simulations presented in this paper, the value of $\varepsilon$ for the Sinkhorn ranking operator does not need to be excessively small. Similarly, for the reflector and refractor costs, FSL demonstrates efficiency gains when $\varepsilon$ is not particularly small. Furthermore, since the size $N$ of the spaces $X$ and $Y$ is typically very large, FSL can significantly reduce the computational complexity of the Sinkhorn iteration in certain scenarios.
\end{remark}

\begin{remark}\label{rmk::multimonial}
    The coefficients $b_{\zeta,\nu}$ often include the multimonial coefficients, see Section~\ref{sec::app_1d} and Section~\ref{sec::app_2d} for examples. Note that all the multimonial coefficients used to compute $b_{\zeta,\nu}$ only need to be calculated for one time, and they can be obtained recursively as follows:
\begin{align*}
    \binom{k}{k_1\,\,k_2\,\cdots\, k_m} = 
	\left\{
	\begin{aligned}
	&1,&\quad&k_m=k,\\
	&\binom{k}{(k_1-1)\,\,k_2\,\cdots\, (k_m+1)}\cdot \frac{k_m+1}{k_1},&\quad&k_1 \ge 1,\\
    &\cdots&&\\
    &\binom{k}{k_1\,\cdots\,(k_{m-1}-1)\,\, (k_m+1)}\cdot\frac{k_m+1}{k_{m-1}},&\quad&k_{m-1} \ge 1.
	\end{aligned}
	\right.
\end{align*}
Therefore, the calculation of each multinomial coefficient only takes $O(1)$.
\end{remark}

\begin{remark}
    In some situations, $b_{\zeta,\nu}$ may be zero for many subscript $\zeta,\nu$, and therefore need not be taken into account. See Section~\ref{sec::app_1d} and Section~\ref{sec::app_2d} for examples.
\end{remark}

\begin{remark}
    The idea of FSL can be easily generalized to high-dimensional cases with $d>1$. In Section~\ref{sec::app_2d}, we will give some examples to illustrate how to implement FSL when $d=2$. The extension to higher-dimensional cases is in a similar way and is therefore omitted. In general, for the $d$-dimensional case, the FSL algorithm performs each Sinkhorn iteration~\eqref{eq::SinkhornIterations} in $O((ML)^dN)$ time when $(ML)^d \ll N$.\footnote{This inequality is naturally satisfied because $N$ is often significantly large in the $d$-dimensional case.}
    \label{prop::main}
\end{remark}

\section{Application to Sinkhorn Ranking}\label{sec::app_1d}

\subsection{Sinkhorn Ranking}
Ranking is a fundamental and important operation used extensively in various areas, such as machine learning, statistics, and information science\cite{cuturi2019differentiable, xie2020differentiable, qin2010general, taylor2008softrank, swezey2021pirank, blondel2020fast}. 

The ranking operator has a close connection with optimal transport. In fact, it can be derived by solving an optimal transport problem\cite{cuturi2019differentiable}. Here, we review the main results of \cite{cuturi2019differentiable}. Given an input vector $\x=(x_1,x_2,\cdots,x_N)^\top \in \mathbb{R}^N$ with distinct entries, a ranking operation returns the ranks of the entries of $\x$. Specifically, it can be seen as an operator $R: \mathbb{R}^N \longrightarrow \mathbb{R}^N$ that outputs a permutation of $\{1,2,\cdots,N\}$, where a higher rank $R(\x)_i$ indicates that $x_i$ has a greater value. The definition of the ranking operator is provided below:
\begin{definition}\label{def::R}
    For a vector $\x=(x_1,x_2,\cdots,x_N)^\top \in \mathbb{R}^N$ whose entries are all distinct, let $\y=(y_1,y_2,\cdots,y_N)^\top \in \mathbb{R}^N$, where $y_1 < y_2 < \cdots < y_N$, $X=\{x_1,\cdots,x_N\}$, $Y=\{y_1,\cdots,y_N\}$, $\bm{a}=\bm{b}= \bm{1}_N/N$, define $C_{ij}=(y_j-x_i)^2$, $i,j=1,2,\cdots,N$, then for any optimal solution $\bm{\Gamma}_*$ to~\eqref{eq::KP}, the ranking operator $R(\x)$ is defined as
\begin{align}
    R(\x) \triangleq N\bm{\Gamma}_*\left(
    \begin{array}{c}
        1\\
        2\\
        \vdots\\
        N
    \end{array}
    \right).
    \label{eq::usual ranking}
\end{align}
\end{definition}

However, the non-differentiability of $R$ limits its applicability in end-to-end learning\cite{cuturi2019differentiable}. The Sinkhorn ranking operator, introduced in\cite{cuturi2019differentiable}, utilizes the entropy regularized optimal transport problem and the Sinkhorn algorithm to enable the differentiability of $R$, thereby facilitating end-to-end training involving ranking and supporting diverse applications\cite{cuturi2019differentiable, xie2020differentiable, blondel2020fast}. The main idea of Sinkhorn ranking is to use the optimal solution to~\eqref{eq::entropyKP} computed by the Sinkhorn algorithm (Algorithm~\ref{alg::Sinkhorn}), which has the form $\mbox{diag}(\bm{\phi}^{(t)}) \bm{K} \mbox{diag}(\bm{\psi}^{(t)})$, to approximate $\bm{\Gamma}_*$. Its definition is as follows:
\begin{definition}\label{def::Sinkhorn R and S}
    For vectors $\x=(x_1,x_2,\cdots,x_N)^\top \in \mathbb{R}^N$ whose entries are all distinct, and $\y=(y_1,y_2,\cdots,y_N)^\top \in \mathbb{R}^N$, where $y_1 < y_2 < \cdots < y_N$, let $X=\{x_1,\cdots,x_N\}$, $Y=\{y_1,\cdots,y_N\}$, $\bm{a} \in \mathcal{P}(X) \cap \mathbb{R}_{++}^N$, $\bm{b} \in \mathcal{P}(Y)\cap \mathbb{R}_{++}^N$, suppose the cost matrix $\bm{C}$ is given by
    \begin{align*}
        C_{ij} = h(y_j-x_i),\quad i,j=1,2,\cdots,N,
    \end{align*}
    where $h:\mathbb{R} \longrightarrow \mathbb{R}_+$ could be an arbitrary non-negative strictly convex differentiable function, given a regularization parameter $\varepsilon>0$, the Sinkhorn ranking operator $R_{h,\varepsilon}(\x, \bm{a}; \y, \bm{b})$ is defined as
    \begin{align}
    R_{h,\varepsilon}(\x, \bm{a}; \y, \bm{b}) \triangleq N\mbox{diag}(\bm{\phi}) \bm{K} \mbox{diag}(\bm{\psi})\left(
    \begin{array}{c}
        b_1\\
        b_1+b_2\\
        \vdots\\
        \sum_{i=1}^N b_i
    \end{array}
    \right)\oslash \bm{a},
    \label{eq::Sinkhorn ranking}
\end{align}
where the vectors $\bm{\phi}$ and $\bm{\psi}$ are the outputs of Algorithm~\ref{alg::Sinkhorn} with initial vector $\bm{\phi}^{(0)} = \bm{1}_N/N$, and $\bm{K}=e^{-\bm{C}/\varepsilon}$, $\bm{b}=(b_1,\cdots,b_N)^\top$.
\end{definition}

It has been proved in\cite{xie2020differentiable} that the Sinkhorn ranking operator $R_{h,\varepsilon}(\x, \bm{a}; \y, \bm{b})$ is differentiable with respect to $\x$. Such a differentiable ranking operator is useful in information retrieval\cite{qin2010general, taylor2008softrank, liu2009learning}.
\begin{example}\label{eg::Sinkhorn ranking}
Consider $\x=(x_1,x_2,x_3)^\top=(0.3, 1.2, -0.25)^\top$, and $\y=(y_1,y_2,y_3)^\top=(0, 0.5, 1)^\top$. Let $\bm{a}=\bm{b}=(1/3,1/3,1/3)^\top$.

In this case, the optimal solution to~\eqref{eq::KP} is 
\begin{align}
    \bm{\Gamma}_* = \frac{1}{3}\left(
    \begin{array}{ccc}
        0 & 1 & 0\\
        0 & 0 & 1\\
        1 & 0 & 0
    \end{array}
    \right).
    \label{eq::R plan}
\end{align}
According to Definition~\ref{def::R}, the ranking operator
\begin{align*}
    R(\x) = 3\bm{\Gamma}_*\left(
    \begin{array}{c}
        1\\
        2\\
        3
    \end{array}
    \right)=\left(
    \begin{array}{c}
        2\\
        3\\
        1
    \end{array}
    \right).
\end{align*}
It can be seen from $R(\x)$ that $x_2>x_1>x_3$ since $R(\x)_2>R(\x)_1>R(\x)_3$.

As for the Sinkhorn ranking operator, let $h(z)=h_{\rm{sq}}(z)=z^2$, which is a non-negative strictly convex differentiable function, then the cost matrix $\bm{C}$ is defined by $C_{ij}=(y_j-x_i)^2$, $i,j=1,2,\cdots,N$. Given $\varepsilon=0.5$, the optimal solution to~\eqref{eq::entropyKP} computed by the Sinkhorn algorithm is
\begin{align}
    \mbox{diag}(\bm{\phi}) \bm{K} \mbox{diag}(\bm{\psi}) \approx \left(
    \begin{array}{ccc}
        0.108 & 0.151 & 0.074\\
        0.010 & 0.081 & 0.242\\
        0.216 & 0.101 & 0.017
    \end{array}
    \right).
    \label{eq::Sinkhorn R plan}
\end{align}
According to Definition~\ref{def::Sinkhorn R and S}, the Sinkhorn ranking operator
\begin{align*}
    R_{h,\varepsilon}(\x, \bm{a}; \y, \bm{b}) = 3\mbox{diag}(\bm{\phi}) \bm{K} \mbox{diag}(\bm{\psi})\cdot\frac{1}{3}\left(
    \begin{array}{c}
        1\\
        2\\
        3
    \end{array}
    \right)\oslash \frac{1}{3}\left(
    \begin{array}{c}
        1\\
        1\\
        1
    \end{array}
    \right) = \left(
    \begin{array}{c}
        1.900\\
        2.698\\
        1.402
    \end{array}
    \right).
\end{align*}

This example is illustrated in Figure~\ref{fig::RS and Sinkhorn RS}. For $i,j=1,2,3$, the width of the arrows from $x_i$ to $y_j$ in Figure~\ref{fig::RS and Sinkhorn RS}(a) and Figure~\ref{fig::RS and Sinkhorn RS}(b) represents the value of the $(i,j)$-th element of the transport plans given in~\eqref{eq::R plan} and~\eqref{eq::Sinkhorn R plan}, respectively.
\begin{figure}[htbp]
    \begin{minipage}[t]{1\textwidth}
        \centering
        \includegraphics[width=.824\linewidth]{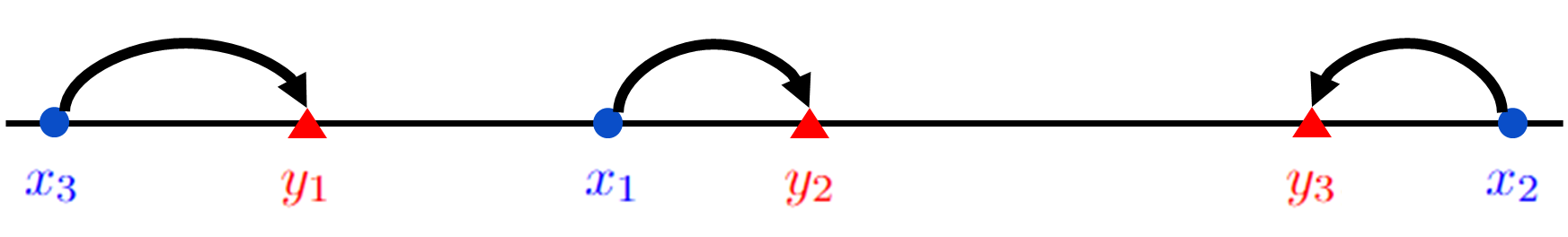}
        \subcaption*{(a)}
    \end{minipage}\\
    
    \begin{minipage}[t]{1\textwidth}
        \centering
        \includegraphics[width=.824\linewidth]{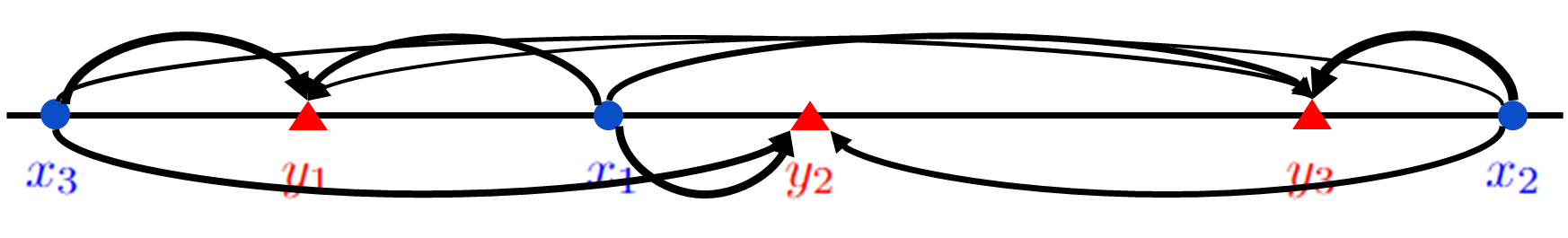}
        \subcaption*{(b)}
    \end{minipage}
    \caption{Illustration of Example~\ref{eg::Sinkhorn ranking}. $\x=(x_1,x_2,x_3)^\top=(0.3, 1.2, -0.25)^\top$, $\y=(y_1,y_2,y_3)^\top=(0, 0.5, 1)^\top$. (a) The ranking operator $R(\x)=(2,3,1)^\top$. It is computed through the transport plan given in~\eqref{eq::R plan}, where $x_i$ is transported to $y_{R(\x)_i}$, $i=1,2,3$. (b) Let $\bm{a}=\bm{b}=(1/3,1/3,1/3)^\top$, $h(z)=h_{\rm{sq}}(z)=z^2$, $\varepsilon=0.5$. The Sinkhorn ranking operator $R_{h,\varepsilon}(\x, \bm{a}; \y, \bm{b})=(1.900,2.698,1.402)^\top$, which is computed through the transport plan given in~\eqref{eq::Sinkhorn R plan}. The width of the arrow from $x_i$ to $y_j$ corresponds to the value of the $(i,j)$-th element of the transport plan in~\eqref{eq::Sinkhorn R plan}, $i,j=1,2,3$.}
    \label{fig::RS and Sinkhorn RS}
\end{figure}
\end{example}

\subsection{Log-type Cost Function and Detailed Implementation of FSL}\label{subsec::1d_FS4Implementation}
While the authors of\cite{cuturi2019differentiable} use $h(z)=z^2$ to define the cost matrix, in this work we assume that $\max_{i=1}^N \{x_i\} < y_1$ and consider the log-type cost function 
\begin{align*}
    h_{\rm{log}}(z) = -\log\left(1-\frac{z}{\tau}\right),\quad 0 < z < \tau,
\end{align*}
where $\tau>y_N-\min_{i=1}^N\{x_i\}$. It can be verified that $h_{\rm{log}}$ is a non-negative strictly convex differentiable function, and therefore, it can be used to define the Sinkhorn ranking operator according to Definition~\ref{def::Sinkhorn R and S}. Additionally, It can be verified that
\begin{align}
    1-\frac{y_j-x_i}{\tau}\in (0,1), \quad\forall i,j=1,2,\cdots,N,
\end{align}
hence $\bm{C}$ is a log-type transport cost associated with the polynomial
\begin{align}
    P(x,y) = 1-\frac{y-x}{\tau},
    \label{eq::P_1d}
\end{align}
as defined in Definition~\ref{def::log-type cost}. This corresponds to the case where $M=1$, and the computation of the Sinkhorn ranking operator can be greatly accelerated by FSL.

In Section~\ref{subsec::1dexperiments}, we will compare these two cost functions and show that our proposed cost function can give better performance in some situations.

In\cite{cuturi2019differentiable}, the authors suggest rescaling the entries of $\x$ to make them fall in $[0,1]$, and setting $\y$ to be the uniform grid on $[0,1]$ before computing the cost matrix $\bm{C}$. Following this idea, we transform $\x$ by
\begin{align}
    \x \mapsto g\left(\frac{\x-\frac{1}{N}(\x^\top \bm{1}_N)\bm{1}_N}{\frac{1}{\sqrt{N}}\|\x-\frac{1}{N}(\x^\top \bm{1}_N)\bm{1}_N\|_2}\right),
    \label{eq::preprocess_x}
\end{align}
where $g$ is the standard logistic function 
\begin{align*}
    g(z)=\frac{1}{1+e^{-z}}.
\end{align*}
We set $\y$ to be the uniform grid on $[1,2]$, namely $y_i = 1+(i-1)/(N-1)$, $i=1,\cdots,N$,
and choose $\tau>y_N-\min_{i=1}^N\{x_i\}$. 

\begin{remark}
The value of $\varepsilon$ does not need to be excessively small to guarantee the accuracy of the Sinkhorn ranking operator. Furthermore, when $\varepsilon$ is not too small, it also ensures the differentiability of the Sinkhorn ranking operator, thereby facilitating its use in training processes that incorporate ranking\cite{cuturi2019differentiable}.
\end{remark}

In what follows, we describe how to implement the FSL algorithm for the log-type transport cost associated with polynomial~\eqref{eq::P_1d}. Suppose that the regularization parameter $\varepsilon=1/L$, where $L$ is a positive integer. In this case, the kernel matrix is given by
\begin{multline*}
    K_{ij} = \left(1-\frac{y_j-x_i}{\tau}\right)^L\\= \sum_{\substack{\zeta,\nu \ge 0\\\zeta+\nu\le L}} \binom{L}{\zeta\,\,\nu\,\,L-\zeta-\nu}\left(\frac{x_i}{\tau}\right)^{\zeta}\left(\frac{-y_j}{\tau}\right)^{\nu}, \quad i,j=1,2,\cdots,N.
\end{multline*}
We denote by
\begin{align*}
    \Tilde{x}_i=\frac{x_i}{\tau},\quad \Tilde{y}_j = \frac{y_j}{\tau}.
\end{align*}
For simplicity, the notations $\Tilde{x}_i$ and $\Tilde{y}_j$ will be abbreviated as $x_i$ and $y_j$, respectively. Then, each element of the kernel can be written as
\begin{align*}
    K_{ij} = \sum_{\zeta=0}^{L}\sum_{\nu=0}^{L} b_{\zeta,\nu}x_i^{\zeta}y_j^{\nu},\quad i,j=1,2,\cdots,N,
\end{align*}
where the the coefficients are
\begin{align*}
    b_{\zeta,\nu} = 
	\left\{
	\begin{aligned}
	&(-1)^{\nu}\binom{L}{\zeta\,\,\nu\,\,L-\zeta-\nu},&\quad&\zeta,\nu \ge 0,\,\, \zeta+\nu \le L,\\
	&0,&\quad&\text{otherwise}.
	\end{aligned}
	\right.
\end{align*}
It should be noted that about half of the coefficients $b_{\zeta,\nu}$ are zero, and we will not compute them. Multimonial coefficients are involved in those non-zero $b_{\zeta,\nu}$, see Remark~\ref{rmk::multimonial} for their calculation. We can then directly apply Algorithm~\ref{alg::FSL_1d} to compute the vectors $\bm{\phi}$ and $\bm{\psi}$ that required in Definition~\ref{def::Sinkhorn R and S}, and obtain the Sinkhorn ranking operator.

\subsection{Numerical Simulations}\label{subsec::1dexperiments}
Now we present a series of numerical simulations to illustrate the effectiveness and efficiency of our approach. We first compare the Sinkhorn ranking operators that are computed by the Sinkhorn algorithm with our proposed log-type cost function $h_{\rm{log}}(z) = -\log(1-z/\tau)$ and the squared distance cost function $h_{\rm{sq}}(z)=z^2$ used in\cite{cuturi2019differentiable}. Then, we contrast the performance of the FSL algorithm and the Sinkhorn algorithm for the log-type transport cost
\begin{align}
    C_{ij} = -\log\left(1-\frac{y_j-x_i}{\tau}\right), \quad i,j=1,2,\cdots,N.
    \label{eq::1d_log-type transport cost}
\end{align}

Consider $\x=(x_1,x_2,\cdots,x_N)^\top=(\theta(1), \theta(2), \cdots, \theta(N))^\top \in \mathbb{R}^N$, where $\theta:\{1,2,\cdots,N\}\rightarrow\{1,2,\cdots,N\}$ is a bijection. As we mentioned in Section~\ref{subsec::1d_FS4Implementation}, we rescale the entries of $\x$ by~\eqref{eq::preprocess_x} before running any algorithm. For the log-type cost and the squared distance cost, we set $y_i = 1+(i-1)/(N-1)$ and $y_i = (i-1)/(N-1)$, $i=1,\cdots,N$, respectively. For the log-type cost function, we set $\tau=(2-\min_{i=1}^N\{x_i\})/(1-1/e)$ to guarantee that $\tau>y_N-\min_{i=1}^N\{x_i\}$ and $C_{ij} \in [0,1]$. 

For the comparison of cost functions $h$, we record the average computational time for the Sinkhorn algorithm to run $1000$ iterations for our proposed log-type cost function $h_{\rm{log}}$ and the squared distance cost function $h_{\rm{sq}}$. After obtaining the Sinkhorn ranking operators $R_{h_{\rm{log}}, \varepsilon}$, $R_{h_{\rm{sq}}, \varepsilon}$ defined in~\eqref{eq::Sinkhorn ranking} and the ranking operator $R$ defined in~\eqref{eq::usual ranking}, say $R_{h_{\rm{log}}, \varepsilon}=(r_1^{\rm{log}},r_2^{\rm{log}},\cdots,r_N^{\rm{log}})^\top$, $R_{h_{\rm{sq}}, \varepsilon}=(r_1^{\rm{sq}},r_2^{\rm{sq}},\cdots,r_N^{\rm{sq}})^\top$, and $R=(r_1,r_2,\cdots,r_N)^\top$, we apply the min-max normalization to scale the range of these operators in $[0,1]$. Specifically, we get $\Tilde{R}_{h_{\rm{log}}, \varepsilon}=(\Tilde{r}_1^{\rm{log}},\cdots,\Tilde{r}_N^{\rm{log}})^\top$, $\Tilde{R}_{h_{\rm{sq}}, \varepsilon}=(\Tilde{r}_1^{\rm{sq}},\cdots,\Tilde{r}_N^{\rm{sq}})^\top$, and $\Tilde{R}=(\Tilde{r}_1,\cdots,\Tilde{r}_N)^\top$, where
\begin{align*}
    \Tilde{r}_i^{\rm{log}} &= \frac{r_i^{\rm{log}}-\min_{j=1}^N\{r_j^{\rm{log}}\}}{\max_{j=1}^N\{r_j^{\rm{log}}\}-\min_{j=1}^N\{r_j^{\rm{log}}\}}, \quad
    \Tilde{r}_i^{\rm{sq}} = \frac{r_i^{\rm{sq}}-\min_{j=1}^N\{r_j^{\rm{sq}}\}}{\max_{j=1}^N\{r_j^{\rm{sq}}\}-\min_{j=1}^N\{r_j^{\rm{sq}}\}},\\
    \Tilde{r}_i &= \frac{r_i-\min_{j=1}^N\{r_j\}}{\max_{j=1}^N\{r_j\}-\min_{j=1}^N\{r_j\}},\quad i=1,2,\cdots,N.
\end{align*}
Then, we compute the mean-squared error (MSE) between $\Tilde{R}_{h, \varepsilon}$ and $\Tilde{R}$ for both $h=h_{\rm{log}}$ and $h=h_{\rm{sq}}$:
\begin{align*}
    \text{MSE}\left(\Tilde{R}_{h_{\rm{log}}, \varepsilon}, \Tilde{R}\right) = \frac{1}{N}\sum_{i=1}^N\left(\Tilde{r}_i^{\rm{log}}-\Tilde{r}_i\right)^2,\quad 
    \text{MSE}\left(\Tilde{R}_{h_{\rm{sq}}, \varepsilon}, \Tilde{R}\right) = \frac{1}{N}\sum_{i=1}^N\left(\Tilde{r}_i^{\rm{sq}}-\Tilde{r}_i\right)^2.
\end{align*}
The length $N$ of the vector $\x$ is set to $200, 400, 800, 1600, 3200, 6400$, and $L$ is set to $10$.

For the comparison of the FSL algorithm and the Sinkhorn algorithm, we conduct two types of tests. On the one hand, we use the FSL algorithm for the log-type transport cost under the same setting as the previous cost function comparison test, and compare the transport plan given by FSL to the one computed by the Sinkhorn algorithm for the log-type transport cost. On the other hand, we record the average computational time of each algorithm for attaining certain marginal error $\|\mbox{diag}(\bm{\psi}^{(t)})\bm{K}^\top \bm{\phi}^{(t)}-\bm{b}\|_1$, where $N$ is set to $10000$, and $L$ is set to $10,15,20$. 

For each algorithm on each test case, we set $\bm{a}=\bm{b}= \bm{1}_N/N$, start from $\bm{\phi}^{(0)}=\bm{1}_N/N$, and perform $30$ independent runs to obtain the average computational time.

\begin{table}[htbp]
\renewcommand{\arraystretch}{1.4}
\begin{center}
\begin{tabular}{c c c c c c}
\hline
\multirow{2}{*}{$N$} & \multicolumn{2}{c}{Computational Time (s)} & & \multicolumn{2}{c}{MSE}\\
\cline{2-3}\cline{5-6}
 &  Log-type & Squared  & & Log-type & Squared \\
\hline
$200$ & \quad $7.30\times 10^{-2}$ & $7.20\times 10^{-2}$ & &$3.46\times 10^{-3}$ &$4.17\times 10^{-3}$ \\
$400$ & \quad $3.18\times 10^{-1}$ & $3.14\times 10^{-1}$ & &$3.49\times 10^{-3}$ & $4.20\times 10^{-3}$\\
$800$ & \quad  $1.57\times 10^{0}$ & $1.55\times 10^{0}$ & &$3.50\times 10^{-3} $& $4.22\times 10^{-3}$\\
$1600$ & \quad $2.01\times 10^{1}$ & $2.00\times 10^{1}$& &$3.51\times 10^{-3} $& $4.23\times 10^{-3}$\\
$3200$ & \quad $9.95\times 10^{1}$ & $9.92\times 10^{1} $& &$3.51\times 10^{-3} $& $4.23\times 10^{-3}$\\
$6400$ & \quad $6.45\times 10^{2}$ & $6.45\times 10^{2} $& &$3.51\times 10^{-3} $& $4.23\times 10^{-3}$\\
\hline
\end{tabular}
\caption{Comparison of the Sinkhorn ranking operators that are computed by the Sinkhorn algorithm with the log-type cost function $h_{\rm{log}}$ and the squared distance cost function $h_{\rm{sq}}$. Column 2-3: for either function $h$, the average computational time of the Sinkhorn algorithm to run $1000$ iterations with $L=10$. Column 4-5: for either function $h$, the MSE between the rescaled $R_{h,\varepsilon}$ and $R$.}
\label{tab::1d_CostComparison}
\end{center}
\end{table}

The results in Table~\ref{tab::1d_CostComparison} demonstrate the competitive performance of the log-type cost function $h_{\rm{log}}$. While the average computational time of the Sinkhorn algorithm using both cost function is comparable, it can be seen that the Sinkhorn ranking operators defined by $h_{\rm{log}}$ consistently yield significantly smaller MSE values compared to the ones defined by $h_{\rm{sq}}$ across all cases. These results indicate that the log-type transport cost~\eqref{eq::1d_log-type transport cost} is potentially a better choice to define Sinkhorn ranking operators in some scenarios. 

\begin{figure}[htbp]
    \centering
    {%
    \includegraphics[width=.47\linewidth]{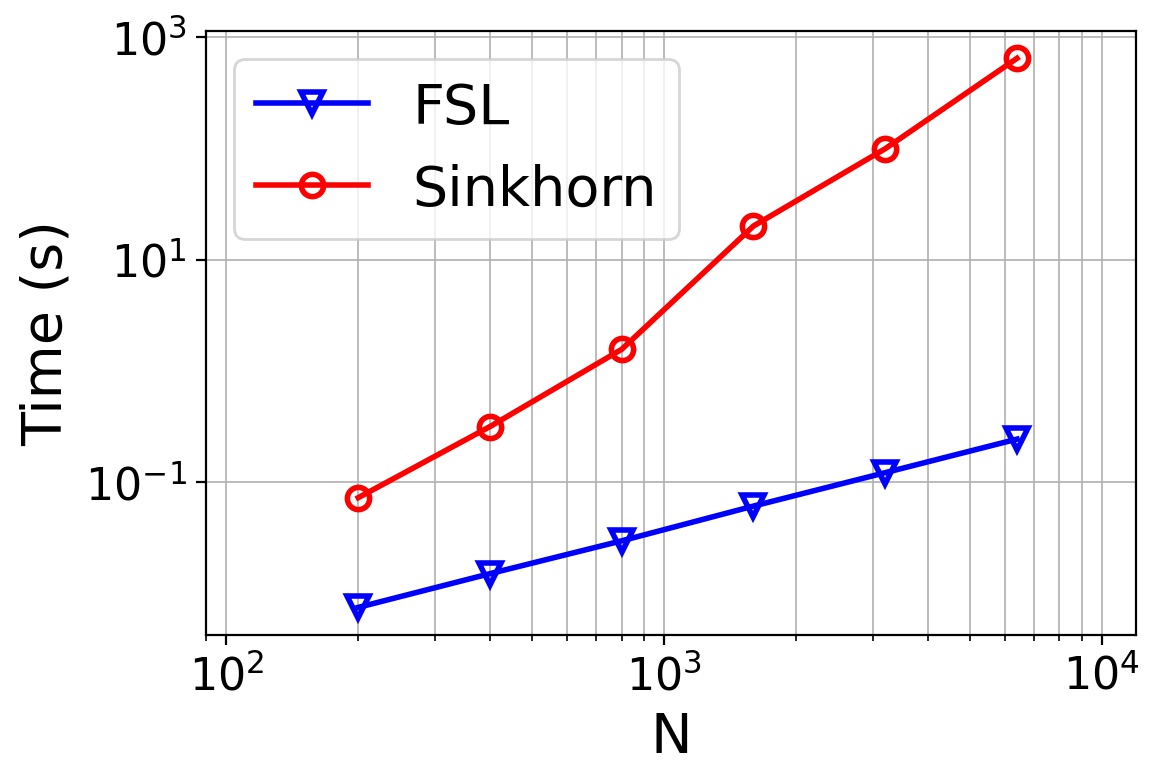}}
    \quad
    {%
    \includegraphics[width=.47\linewidth]{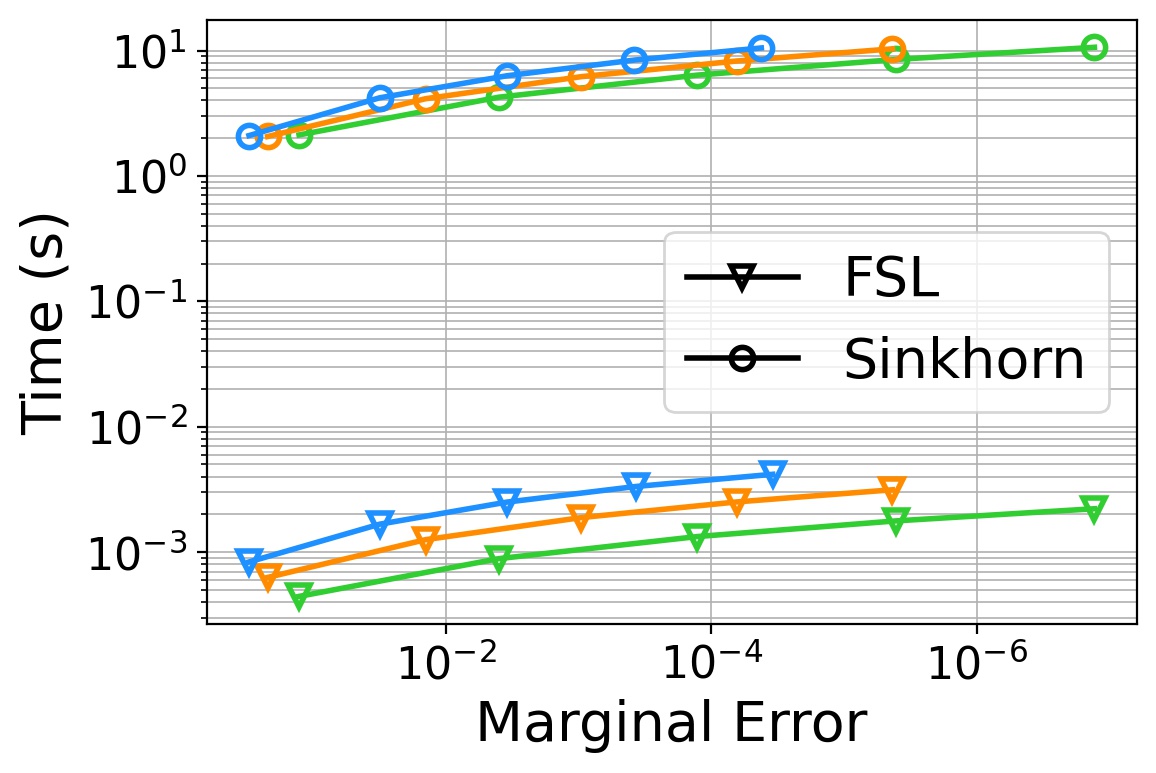}}
    \caption{Comparison of the FSL algorithm and the Sinkhorn algorithm for the log-type cost function. Left: The average computational time of these two algorithms for performing $1000$ iterations, with the same setting as in Table~\ref{tab::1d_TimeComparison}. Right: The average computational time of both algorithms for reaching the corresponding marginal error with $L = 10$ (Green), $L = 15$ (Orange), $L = 20$ (Blue), and $N=10000$.}
    \label{fig::1d_softsorting}
\end{figure}

The comparison between the FSL algorithm and the Sinkhorn algorithm is given in Table~\ref{tab::1d_TimeComparison} and Figure~\ref{fig::1d_softsorting}. From Table~\ref{tab::1d_TimeComparison} and Figure~\ref{fig::1d_softsorting} (Left), it can be found that the FSL algorithm is faster than the Sinkhorn algorithm for all given values of $N$, which shows the high efficiency of the FSL algorithm. The computational time of the FSL algorithm scales approximately linear with $N$, while the one of the Sinkhorn algorithm grows much more rapidly as $N$ increases. The last column of Table~\ref{tab::1d_TimeComparison} shows that the transport plans obtained by the two algorithms are almost identical in all cases. In Figure~\ref{fig::1d_softsorting} (Right), we present the average computational time of both algorithms for reaching the corresponding marginal error. It is evident that the FSL algorithm is more efficient than the Sinkhorn algorithm for different values of $L$.

In summary, the log-type cost function $h_{\rm{log}}$ is sometimes a favorable choice for the Sinkhorn ranking operator, and the computation of this operator can be substantially accelerated by the FSL algorithm. Consequently, our approach offers promising application utility for differentiable ranking.

\begin{table}[htbp]
\renewcommand{\arraystretch}{1.4}
\begin{center}
\begin{tabular}{c c c c c}
\hline
\multirow{2}{*}{$N$} & \multicolumn{2}{c}{Computational Time (s)} & \multirow{2}{*}{Speed-up Ratio} & \multirow{2}{*}{$\|\bm{\Gamma}_{\rm{FSL}}-\bm{\Gamma}_{\rm{S}}\|_F$}\\
\cline{2-3}
 & \quad FSL & Sinkhorn & & \\
\hline
$200$ & \quad $7.55\times 10^{-3}$ & $7.30\times 10^{-2}$ & $9.67\times 10^{0}$ &$8.99\times 10^{-14}$ \\
$400$ & \quad $1.52\times 10^{-2}$ & $3.18\times 10^{-1}$ & $2.09\times 10^{1}$ & $2.01\times 10^{-14}$\\
$800$ & \quad  $2.99\times 10^{-2}$ & $1.57\times 10^{0}$ & $5.25\times 10^{1} $& $1.03\times 10^{-14}$\\
$1600$ & \quad $6.14\times 10^{-2}$ & $2.01\times 10^{1} $& $3.28\times 10^{2} $& $3.90\times 10^{-15}$\\
$3200$ & \quad $1.22\times 10^{-1}$ & $9.95\times 10^{1} $& $8.13\times 10^{2} $& $7.34\times 10^{-15}$\\
$6400$ & \quad $2.45\times 10^{-1}$ & $6.45\times 10^{2} $& $2.63\times 10^{3} $& $7.13\times 10^{-15}$\\
\hline
\end{tabular}
\caption{Comparison of the FSL algorithm and the Sinkhorn algorithm for the log-type cost function. Column 2-3: the average computational time of these two algorithms to execute $1000$ iterations with $L=10$. Column 4: the ratio of the average computational time of the Sinkhorn algorithm to that of the FSL algorithm. Column 5: the Frobenius norm of the difference of the transport plans obtained by the two algorithms, $\bm{\Gamma}_{\rm{FSL}}$ represents the transport plan obtained by the FSL algorithm, and $\bm{\Gamma}_{\rm{S}}$ represents the one obtained by the Sinkhorn algorithm.}
\label{tab::1d_TimeComparison}
\end{center}
\end{table}

\section{Application to Reflector and Refractor Costs}\label{sec::app_2d}

\subsection{Reflector and Refractor Costs}
We consider in what follows $X,Y \subset \mathbb{R}^2$. The reflector cost is the log-type transport cost derived from the far-field reflector problem\cite{wang2004design}, which is given by
\begin{align*}
    C_{ij} = -\log(1-\langle \x_i,\y_j\rangle),\quad i,j=1,2,\cdots,N.
\end{align*}
The refractor cost is the transport cost that derived from the far-field refractor problem with the relative refractive index $0<\kappa < 1$\cite{gutierrez2009refractor}, which is written as
\begin{align*}
    C_{ij} = -\log(1-\kappa\langle \x_i,\y_j\rangle),\quad i,j=1,2,\cdots,N,
\end{align*}
where $\kappa$ is the ratio of the index of refraction of the media that contains the target region to the one of the media that contains the point source, as introduced in\cite{gutierrez2009refractor}. 

For convenience, we write the transport cost in both cases as
\begin{align}
    C_{ij} = -\log(1-\kappa\langle \x_i,\y_j\rangle),\quad i,j=1,2,\cdots,N,
    \label{eq::R-type}
\end{align}
where $\kappa=1$ represents the reflector cost, and $0<\kappa < 1$ represents the refractor cost.

\subsection{Detailed Implementation of FSL}
The corresponding polynomial of the reflector and refractor cost~\eqref{eq::R-type} is taken as
\begin{align*}
    P(x_{i,1},x_{i,2},y_{j,1},y_{j,2}) = 1-\kappa x_{i,1}y_{j,1}-\kappa x_{i,2}y_{j,2}, \quad i,j=1,2,\cdots,N,
\end{align*}
which is in the form of~\eqref{eq::P} with $M=1$.
Assume that~\eqref{eq::P_assume} holds and the regularization parameter $\varepsilon = 1/L$ with $L \in \mathbb{N}$, the computational cost for the FSL algorithm to compute $\bm{\phi}^{(t+1)}$ and $\bm{\psi}^{(t+1)}$ is therefore $O(L^2N)$ when $L^2 \ll N$ according to Remark~\ref{prop::main}.\footnote{In the two-dimensional case, grid-like sampling is commonly employed, with a substantial average sampling size in each dimension, denoted by $N_1$ and $N_2$. Therefore, we naturally have $L \ll N_1$ and $L \ll N_2$, and thus $L^2 \ll N_1 N_2 \approx N$.}

Below, we explain how to apply the FSL algorithm for this transport cost. In this situation, each element of the kernel $\bm{K}$ is written as
\begin{multline*}
    K_{ij} = \left(1-\kappa x_{i,1}y_{j,1}-\kappa x_{i,2}y_{j,2}\right)^L\\
    =\sum_{\substack{p,q \ge 0\\p+q\le L}} (-\kappa)^{p+q}\binom{L}{p\,\,q\,\,L-p-q} (x_{i,1}y_{j,1})^p(x_{i,2}y_{j,2})^q,
    \quad i,j=1,2,\cdots,N.
\end{multline*}
Let
\begin{align*}
    b_{pq} \triangleq (-\kappa)^{p+q}\binom{L}{p\,\,q\,\,L-p-q},\quad  p,q \ge 0,\,\,p+q \le L.
\end{align*}
For $z=x$ or $y$, we denote by
\begin{align*}
    \bm{A}^z &= (A^z_{pq,j})\triangleq \left(
    \begin{array}{cccc}
        1 & 1 & \cdots & 1\\
        z_{1,2} & z_{2,2} & \cdots & z_{N,2}\\
        z_{1,1} & z_{2,1} & \cdots & z_{N,1}\\
        \vdots & \vdots & \ddots & \vdots\\
        z_{1,2}^L & z_{2,2}^L & \cdots & z_{N,2}^L\\
        z_{1,1}z_{1,2}^{L-1} & z_{2,1}z_{2,2}^{L-1} & \cdots & z_{N,1}z_{N,2}^{L-1}\\
        \vdots & \vdots & \ddots & \vdots\\
        z_{1,1}^L & z_{2,1}^L & \cdots & z_{N,1}^L\\
    \end{array}
    \right)\in \mathbb{R}^{\frac{(L+1)(L+2)}{2}\times N},\nonumber\\
    A^z_{pq,j} &\triangleq z_{j,1}^pz_{j,2}^q,\quad p,q \ge 0,\,\,p+q \le L,\,\, j=1,2,\cdots,N.
\end{align*}
Then, each entry of $\bm{K}\bm{\xi}$ becomes
\begin{align}
    (\bm{K} \bm{\xi})_i = \sum_{\substack{p,q \ge 0\\p+q\le L}}b_{pq}A^x_{pq,i}\left(\bm{A}^y\bm{\xi}\right)_{pq},\quad i=1,2,\cdots,N.
    \label{eq::FSR}
\end{align}
Similar to Algorithm~\ref{alg::FSL_1d}, FSL first computes $\bm{A}^y\bm{\xi}$, and then obtains each entry of $\bm{K}\bm{\xi}$ by~\eqref{eq::FSR}.

Note that the above procedure is simpler than that in Section~\ref{subsec::FSL}. This is because many coeffeients $b_{\zeta_1 ,\zeta_2 ,\nu_1 ,\nu_2 }$ are equal to zero in the expansion of $K_{ij}$. Indeed, $b_{\zeta_1 ,\zeta_2 ,\nu_1 ,\nu_2 } \ne 0$ only when $p=\zeta_1 =\nu_1 $, $q=\zeta_2 =\nu_2 $, and $p+q \le L$. 

\subsection{Numerical Simulations}
Let $X=Y \subset (0,1)^2$ be a two-dimensional uniform square grid of size $N=n \times n$ with grid length $h=0.7/(1.1n)$. The grid points are $(ih, jh)$, $i,j=1,2,\cdots,n$. We use the FSL algorithm and the Sinkhorn algorithm to solve the entropy regularized problem~\eqref{eq::entropyKP} between two random probability measures $\bm{a},\bm{b} \in \mathcal{P}(X)$, where the probability vectors are constructed as follows:
\begin{enumerate}
    \item Let $\bm{a}=(a_1,a_2,\cdots,a_N)^\top$, $\bm{b}=(b_1,b_2,\cdots,b_N)^\top$ with $a_i,b_i$ i.i.d.$\sim \mbox{U}(0,1)$, $i=1,2,\cdots,N$.
    \item Set $\bm{a} \leftarrow \frac{\bm{a}}{\|\bm{a}\|_1}$, $\bm{b} \leftarrow \frac{\bm{b}}{\|\bm{b}\|_1}$.
\end{enumerate}
The transport costs are the reflector cost and the refractor cost given by~\eqref{eq::R-type}. For the refractor case, we set $\kappa=0.5$. 

We use almost the same settings as in Section~\ref{subsec::1dexperiments}, except that the number $N$ of grid points is set to $20\times 20, 40\times 40, 80\times 80, 160\times 160$ for the $1000$-iteration test, and is set to $100 \times 100$ for the marginal error test. 

\begin{figure}[htbp]
    \begin{minipage}[t]{1\textwidth}
        \centering
        \subcaption*{(a) Reflector Cost}
        {%
        \includegraphics[width=.47\linewidth]{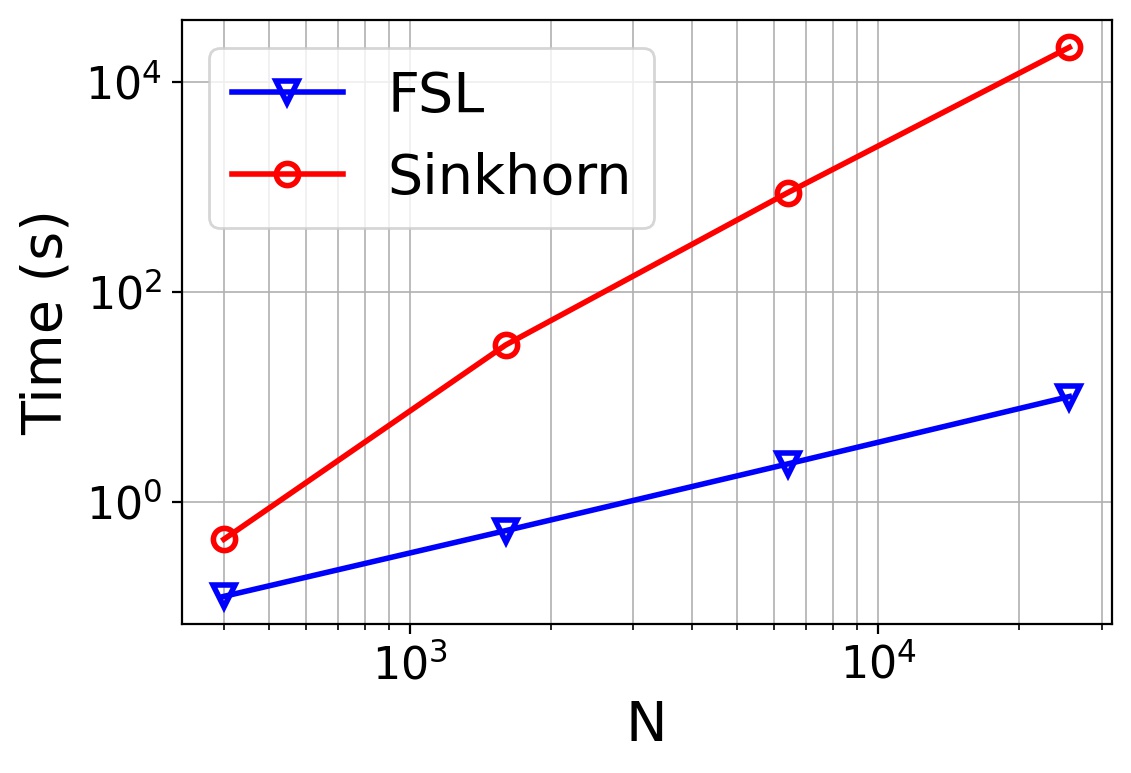}}
        \quad
        {%
        \includegraphics[width=.47\linewidth]{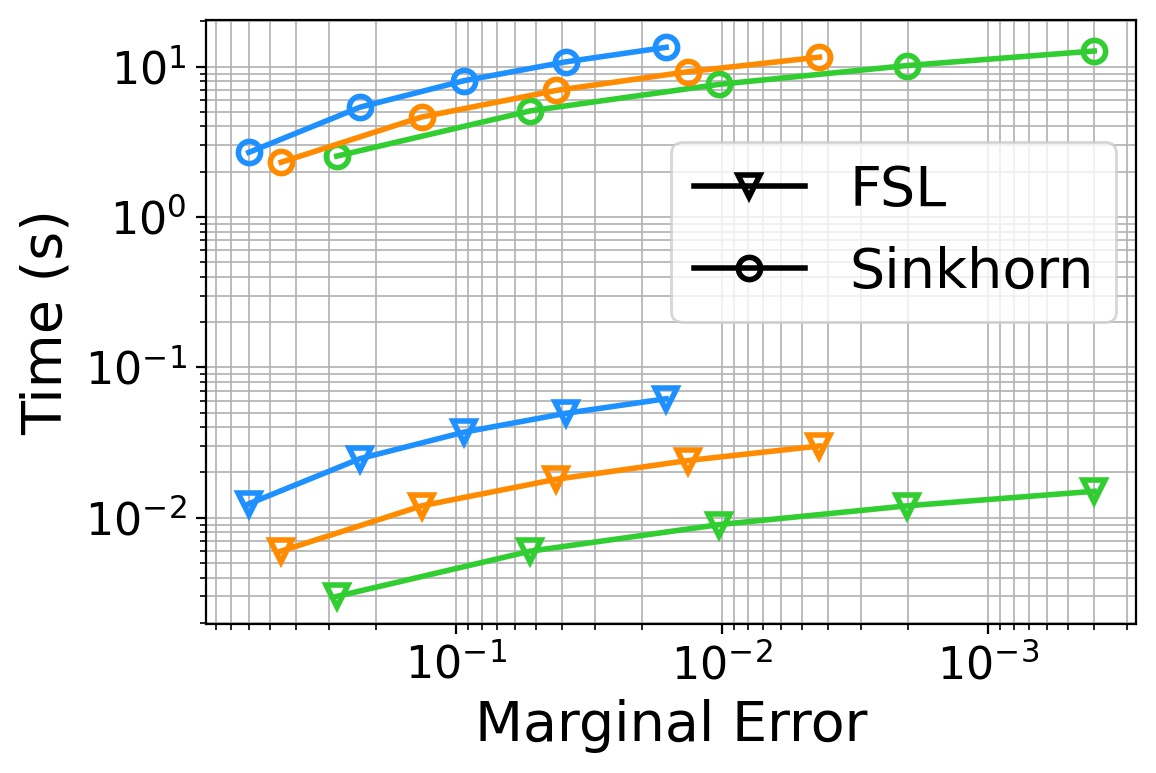}}
    \end{minipage}\\
    \begin{minipage}[t]{1\textwidth}
        \centering
        \subcaption*{(b) Refractor Cost ($\kappa=0.5)$}
        {%
        \includegraphics[width=.47\linewidth]{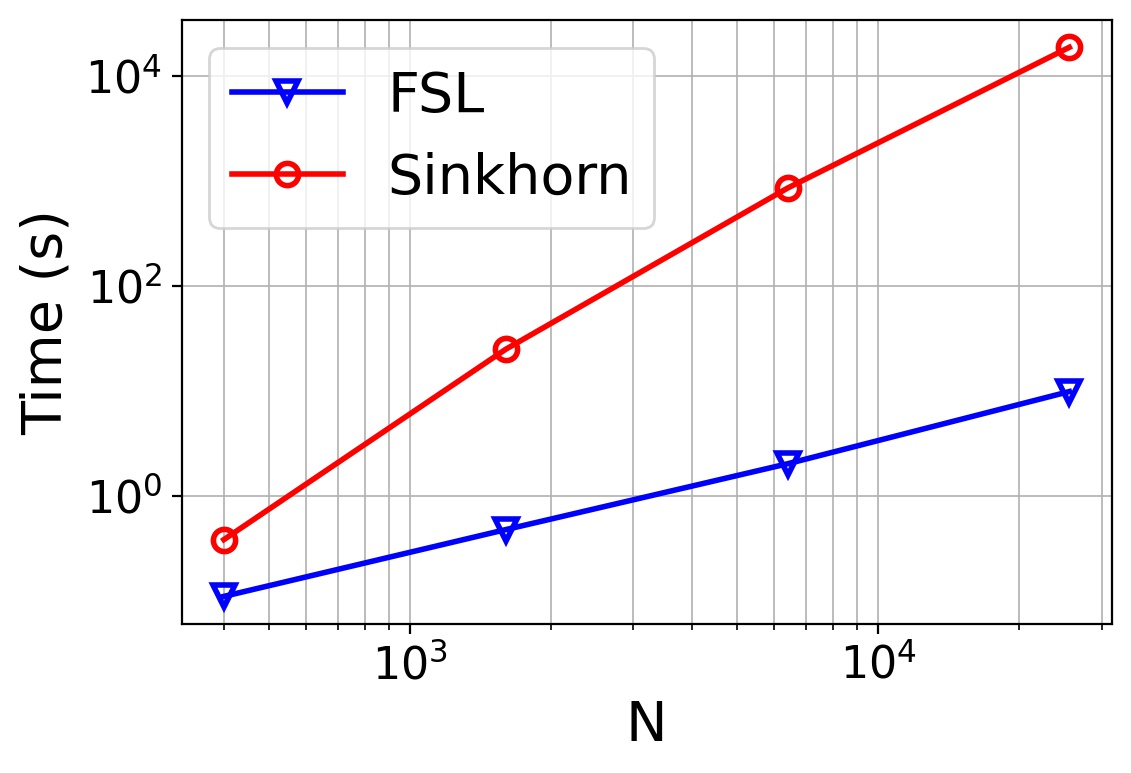}}
        \quad
        {%
        \includegraphics[width=.47\linewidth]{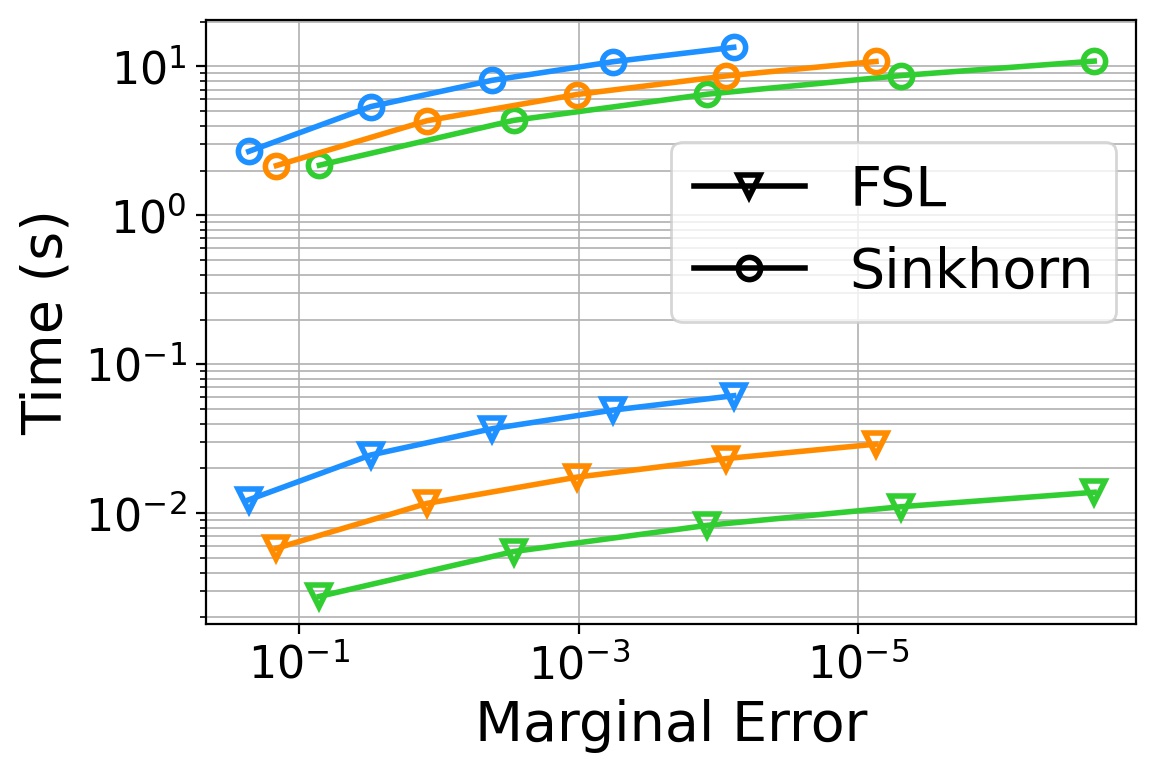}}
    \end{minipage}
    \caption{Comparison of the FSL algorithm and the Sinkhorn algorithm for the reflector and refractor costs. (a) Test for the reflector cost. (b) Test for the refractor cost ($\kappa=0.5$). Left: The average computational time for performing $1000$ iterations. Right: The average computational time for reaching the corresponding marginal error with $L = 10$ (Green), $L = 15$ (Orange), $L = 20$ (Blue), and $N=100\times 100$ grid points.}
    \label{fig::2d}
\end{figure}

The results for the reflector cost are presented in Table~\ref{tab::2d_reflector} and Figure~\ref{fig::2d} (a), and the results for the refractor cost are given in Table~\ref{tab::2d_refractor} and Figure~\ref{fig::2d} (b). These results are similar to that in Section~\ref{subsec::1dexperiments}. It can be observed that, in all cases, $\Gamma_{FSL}$ is virtually the same as $\Gamma_{S}$, and the FSL algorithm needs far less time than the Sinkhorn algorithm. In summary, the FSL algorithm shows substantial advantages in terms of efficiency in comparison to the Sinkhorn algorithm for both the reflector and refractor costs. 

\begin{table}[htbp]
\renewcommand{\arraystretch}{1.4}
\begin{center}
\begin{tabular}{c c c c c}
\hline
\multirow{2}{*}{$N=n \times n$} & \multicolumn{2}{c}{Computational Time (s)} & \multirow{2}{*}{Speed-up Ratio} & \multirow{2}{*}{$\|\bm{\Gamma}_{\rm{FSL}}-\bm{\Gamma}_{\rm{S}}\|_F$}\\
\cline{2-3}
 & \quad FSL & Sinkhorn & & \\
\hline
$20 \times 20$ & \quad $1.25\times 10^{-1}$ & $4.37\times 10^{-1}$ & $3.50\times 10^{0}$ & $1.30\times 10^{-15}$\\
$40 \times 40$ & \quad  $5.28\times 10^{-1}$ & $3.11\times 10^{1}$ & $5.89\times 10^{1} $& $6.13\times 10^{-16}$\\
$80 \times 80$ & \quad $2.28\times 10^{0}$ & $8.78\times 10^{2} $& $3.85\times 10^{2} $& $2.64\times 10^{-16}$\\
$160 \times 160$ & \quad $1.00\times 10^{1}$ & $2.14\times 10^{4} $& $2.14\times 10^{3} $& $3.27\times 10^{-17}$\\
\hline
\end{tabular}
\caption{Comparison of the FSL algorithm and the Sinkhorn algorithm for the reflector cost. Column 2-3: the average computational time of these two algorithms to execute $1000$ iterations with $L=10$. Column 4: the ratio of the average computational time of the Sinkhorn algorithm to that of the FSL algorithm. Column 5: the Frobenius norm of the difference of the approximate transport
plans obtained by the two algorithms.}
\label{tab::2d_reflector}
\end{center}
\end{table}
\begin{table}[htbp]
\renewcommand{\arraystretch}{1.4}
\begin{center}
\begin{tabular}{c c c c c}
\hline
\multirow{2}{*}{$N=n \times n$} & \multicolumn{2}{c}{Computational Time (s)} & \multirow{2}{*}{Speed-up Ratio} & \multirow{2}{*}{$\|\bm{\Gamma}_{\rm{FSL}}-\bm{\Gamma}_{\rm{S}}\|_F$}\\
\cline{2-3}
 & \quad FSL & Sinkhorn & & \\
\hline
$20 \times 20$ & \quad $1.10\times 10^{-1}$ & $3.83\times 10^{-1}$ & $3.47\times 10^{0}$ & $4.08\times 10^{-17}$\\
$40 \times 40$ & \quad  $4.76\times 10^{-1}$ & $2.48\times 10^{1}$ & $5.21\times 10^{1} $& $7.76\times 10^{-18}$\\
$80 \times 80$ & \quad $2.02\times 10^{0}$ & $8.51\times 10^{2} $& $4.22\times 10^{2} $& $3.20\times 10^{-18}$\\
$160 \times 160$ & \quad $9.83\times 10^{0}$ & $1.87\times 10^{4} $& $1.90\times 10^{3} $& $1.31\times 10^{-18}$\\
\hline
\end{tabular}
\caption{Comparison of the FSL algorithm and the Sinkhorn algorithm for the refractor cost ($\kappa=0.5$). Column 2-3: the average computational time of these two algorithms to execute $1000$ iterations with $L=10$. Column 4: the ratio of the average computational time of the Sinkhorn algorithm to that of the FSL algorithm. Column 5: the Frobenius norm of the difference of the approximate transport
plans obtained by the two algorithms.}
\label{tab::2d_refractor}
\end{center}
\end{table}

\section{Conclusion}\label{sec::conclusion}
In this paper, we put forth the FSL algorithm, which can achieve linear-time complexity for each Sinkhorn iteration, for the optimal transport problems with a specific kind of log-type transport cost. This algorithm has broad application prospects, as it can be utilized, together with a log-type cost function, to compute the Sinkhorn ranking operator; and can be applied to the reflector and refractor costs. While we have only considered the cases on the plane, the real-world far-field reflector and refractor problems are set on the sphere. The extension of our method to spherical cases is still under investigation, and we hope to report our progress in the future work.

\section*{Acknowledgment}
This work was supported by National Natural Science Foundation of China (Grant No. 12271289). The authors would like to express their gratitude to Prof. Xujia Wang of Westlake University for his valuable discussions and suggestions.

\end{document}